# Towards a Global Controller Design for Guaranteed Synchronization of Switched Chaotic Systems

Indranil Pan[a], Saptarshi Das[a,b*], and Avijit Routh[a]

a) *Department of Power Engineering, Jadavpur University, Salt Lake Campus, LB-8, Sector 3, Kolkata-700098, India.*
b) *School of Electronics and Computer Science, University of Southampton, Southampton SO17 1BJ, United Kingdom.*

**Authors' emails:**
indranil.jj@student.iitd.ac.in, indranil@pe.jusl.ac.in (I. Pan)
saptarshi@pe.jusl.ac.in, s.das@soton.ac.uk (S. Das*)
avijit_hit05@yahoo.co.in (A. Routh)

Phone no. +44(0)7448572598

**Abstract:**
In this paper, synchronization of identical switched chaotic systems is explored based on Lyapunov theory of guaranteed stability. Concepts from robust control principles and switched linear systems are merged together to derive a sufficient condition for synchronization of identical master-slave switched nonlinear chaotic systems and are expressed in the form of bilinear matrix inequalities (BMIs). The nonlinear controller design problem is then recast in the form of linear matrix inequalities (LMIs) to facilitate numerical computation by standard LMI solvers and is illustrated by appropriate examples.

**Keywords:** chaos synchronization; Linear Matrix Inequality (LMI); nonlinear state feedback controller; robust stability; switched chaotic system; unified chaotic system

## 1. Introduction

Synchronization of two chaotic systems was first introduced and experimentally demonstrated by Carroll and Pecora [1]. Since then chaos synchronization based schemes have found numerous applications in secure communication [2], cryptography [3], chaotic shift keying [4], chaotic modulation [5] etc. For a given chaotic system which is termed as the master or driving system, the concept of synchronization is to make the states of another identical chaotic system (termed as slave or response system) follow those of the master system. Many chaos synchronization methods have been investigated by contemporary researchers, like the linear feedback control [6], time delayed feedback control [7], artificial intelligence based optimization methods [8], nonlinear control [9], active control [10], adaptive sliding mode control [11] etc. In this paper, a nonlinear controller design technique is proposed in terms of BMI formulation for synchronization of switched chaotic systems, which are then recast as a set of LMIs. This makes the global nonlinear controller design procedure easy as standard numerical solvers [12], [13] employing efficient interior point algorithms are available for solving such LMI problems in polynomial time.

Hybrid systems are ones in which both the continuous and discrete dynamics co-exist. Switched systems are a class of hybrid systems in which the discrete dynamics of the switching function is neglected [14]. A switching rule dictates the sequence in which the subsystems are switched and the time for which each subsystem is active. This switching rule can be arbitrary where the user has no control over it, or it may be controlled according to the designer's specifications. Hybrid systems are becoming more important due to their ability to model a large class of problems in biology [15], control systems [16], chemical reactions [17], communication networks [18] etc. Theoretical results regarding switched systems can be found in [14], [19], [20].

In this paper, we consider the case of synchronization of a unified chaotic system [21] where the master and slave systems simultaneously switch from one chaotic system to the other. The switching rule is not predefined and is assumed to be arbitrary. The objective is to design a single nonlinear controller so that the synchronization of the master and slave system is maintained in spite of the unknown arbitrary switching amongst various chaotic systems. Very few works have been done on switched chaotic systems [18], and the motivation of this paper is to develop a universal controller which synchronizes all members of the family of unified nonlinear chaotic systems even under arbitrary switching amongst them. Under such arbitrary switching among the unified chaotic system family, the synchronizing hardware i.e. the controller structure does not need to be changed every time. In many applications, the synchronizing control circuit can only be made once and should not be changed with the nature of the chaotic system. Designing dedicated synchronizing controllers for each different chaotic system may be infeasible to implement in many realistic scenario. Herein lies the applicability of designing a generalized global controller which can handle arbitrary switching among different chaotic system. From hardware implementation point of view, it is also economic to have a global controller. Such a controller, once designed can be treated as a plug and play device for synchronization between any chaotic systems amongst the unified chaotic system family.

There have been several efforts to study synchronization between similar and different chaotic systems e.g. using feedback control, time delay methods, fuzzy, impulsive control [22], [23] and robust synchronization [24] for various applications like chaos based secure communication [25] etc. But most of the papers deal with a fixed structure of the chaotic system from which the error dynamical system can easily be obtained analytically. A stabilization scheme for the error dynamical system can then be used to ensure successful master-slave chaos synchronization. On the other hand, there have been huge advancements in the last few years on the application of robust stability analysis for switched linear systems and hybrid systems [14][26][27]. Lyapunov stabilization problem for a set of linear systems with arbitrary switching amongst them can be solved in terms of few set of LMIs which are derived from a successful bridging between the switched linear systems theory and robust stability theory. Unfortunately, robust stabilization problems for switching amongst complex nonlinear systems are not adequately addressed in literature and to the best of our knowledge, such results for switched nonlinear chaotic systems does not exist. In this paper, we formulate a synchronization scheme between two chaotic systems which continuously switches its characteristics within the Unified chaotic system family. Therefore, here we propose a new class of chaotic systems known as switched chaotic systems and also design a global control scheme to ensure guaranteed synchronization with the consideration of arbitrary switching in the master and slave chaotic system. Here, we study few interesting cases of switching phenomena in the chaotic system like switching of the key parameter ($\alpha$) of the Unified chaotic system, along with the switching of the initial conditions, on-off switching of the synchronizing nonlinear controller etc. It is shown that all of these cases can be synchronized with only one value of the controller gains. Hence this might be termed as a 'global' controller.

LMI problems are convex in nature and there are efficient interior point methods [28] to numerically solve such problems. Lyapunov stabilization scheme for switched systems often reduces to bilinear matrix inequality (BMIs) but these problems are generally non-convex in nature which cannot be solved directly using LMI solvers. In addition, the BMIs are mostly 'NP-hard' problems [29] where the class 'P' indicates problems solvable in polynomial time. Reducing the computational complexity to the solution of the NP hard problems is still an active research area. It is known that the NP-hard is not a characteristics for a particular algorithm but of the problem itself and in many cases such problems could be solved using some approximation or heuristics [30]. For example, the branch and bound algorithm [30] can be employed to solve BMIs but the objective function in such cases need to have a tight upper and lower bound. Also, it is well known that many problems encountered in control theory can be cast as BMI problems and there are many ways to solve BMIs. As in the case of LMIs, the computationally efficient interior point methods cannot be directly employed to solve BMIs. But depending on the structure of the control problem formulation, some mathematical simplification could yield a set of LMIs which can be easily solved. However, it is important to emphasize that such simplification of BMIs to be solvable using standard LMI tools using some analytical treatment or mathematical transformation is extremely problem dependent and as such cannot be considered as a generalized framework to solve all BMI problems. Otherwise BMI solution techniques use some sort of local search techniques which depends on the initial guess values and converges to a local minima. In such an approach there is no guarantee that the obtained solution is a global optima [31][32][33] e.g. the D–K iteration for $\mu$-synthesis [34], alternating semi-definite programs (SDP) method [31], dual iteration method [35]. More details on the computational complexity of BMI problems arising in control and stabilization problems of switched systems and the use of coupled heuristic optimization and convex optimization or the memetic algorithms has been illustrated in Pan and Das [36]. However, in the present problem, the BMI formulation has been reduced to a set of LMIs with suitable transformation [37] so that the efficient interior point method based solution techniques can be applied. The robust control toolbox of Matlab [12] can be employed to solve the LMIs. The other popular software package YALMIP [13] gives more flexibility in implementing the LMIs and hence is widely used in the control community. Simulation and optimization based optimal controller design has been previously investigated in Das *et al.* [38] with a particular choice of initial conditions for establishing synchronization amongst two chaotic systems or the suppression of chaotic oscillation in highly complex systems [39], without the consideration of switching. But such synchronization is not theoretically guaranteed for all possible initial conditions. This is addressed in the present paper with the Lyapunov stabilization and a BMI formulation for the switched chaotic system scenario.

The rest of the paper is organized as follows. Section 2 gives an expository overview of the unified chaotic system, synchronization of chaotic systems and switched systems. In Section 3 sufficient conditions for synchronization between switched chaotic systems are derived using robust control principles. The controller design problem is also addressed and expressed in the form of LMIs. Section 4 gives credible numerical simulations illustrating the design methodology. It also looks into the results and discussion of the proposed design methodology. The paper ends in Section 5, followed by the references.

## 2. Basics of the unified chaotic system, its synchronization and switched systems

The unified chaotic system is given by equation (1) [21].

$$
\begin{aligned}
\dot{x} &= (25\alpha + 10)(y - x) \\
\dot{y} &= (28 - 35\alpha)x + (29\alpha - 1)y - xz \\
\dot{z} &= xy - \left(\frac{\alpha + 8}{3}\right)z
\end{aligned}
\tag{1}
$$

where, the key parameter $\alpha \in [0,1]$. The dynamical system given by (1) is chaotic for all values of $\alpha \in [0,1]$. When $\alpha \in [0,0.8)$, the system (1) is called the general family of Lorenz systems. When $\alpha = 0.8$, it is the Lu system. When $\alpha \in (0.8,1]$, the system is the general family of Chen systems.

The drive or the master system (2) is denoted by the subscript '$m$'.

$$
\begin{aligned}
\dot{x}_m &= (25\alpha + 10)(y_m - x_m) \\
\dot{y}_m &= (28 - 35\alpha)x_m + (29\alpha - 1)y_m - x_m z_m \\
\dot{z}_m &= x_m y_m - \left(\frac{\alpha + 8}{3}\right)z_m
\end{aligned}
\tag{2}
$$

The corresponding slave or response system (3) is denoted by the subscript '$s$',

$$
\begin{aligned}
\dot{x}_s &= (25\alpha + 10)(y_s - x_s) + u_1 \\
\dot{y}_s &= (28 - 35\alpha)x_s + (29\alpha - 1)y_s - x_s z_s + u_2 \\
\dot{z}_s &= x_s y_s - \left(\frac{\alpha + 8}{3}\right)z_s + u_3
\end{aligned}
\tag{3}
$$

where, $u_1, u_2, u_3$ are the non linear control inputs to the three states of the slave system.

The synchronization error is given by (4)

$$e_1 = x_s - x_m$$
$$e_2 = y_s - y_m \quad (4)$$
$$e_3 = z_s - z_m$$

and the corresponding error dynamical system is defined by (5).

$$\dot{e}_1 = (25\alpha + 10)(e_2 - e_1) + u_1$$
$$\dot{e}_2 = (28 - 35\alpha)e_1 + (29\alpha - 1)e_2 - z_m e_1 - x_s e_3 + u_2 \quad (5)$$
$$\dot{e}_3 = y_m e_1 + x_s e_2 - \left(\frac{\alpha + 8}{3}\right)e_3 + u_3$$

The state space representation of the above system in (5) is shown in (6).

$$\begin{bmatrix} \dot{e}_1 \\ \dot{e}_2 \\ \dot{e}_3 \end{bmatrix} = \begin{bmatrix} -(25\alpha + 10) & (25\alpha + 10) & 0 \\ (28 - 35\alpha) - z_m & (29\alpha - 1) & -x_s \\ y_m & x_s & -\left(\frac{\alpha + 8}{3}\right) \end{bmatrix} \begin{bmatrix} e_1 \\ e_2 \\ e_3 \end{bmatrix} + \begin{bmatrix} u_1 \\ u_2 \\ u_3 \end{bmatrix}$$

(6)

$$\Rightarrow \dot{e} = Ae + f(u)$$

where, $A = \begin{bmatrix} -(25\alpha + 10) & (25\alpha + 10) & 0 \\ (28 - 35\alpha) - z_m & (29\alpha - 1) & -x_s \\ y_m & x_s & -\left(\frac{\alpha + 8}{3}\right) \end{bmatrix}$ and $f(u) = \begin{bmatrix} u_1 \\ u_2 \\ u_3 \end{bmatrix}$

Since the error dynamical system (6) is defined as the difference of the slave and master states, the asymptotic stability of the nonlinear system (6) implies that the master and the slave states are synchronized.

In mathematical terms, a switched system can be defined by a differential equation [20] as shown in (7).

$$\dot{x} = g_\sigma(x) \quad (7)$$

where, $\left\{ g_p : p \in \Lambda \right\}$ is a family of sufficiently regular functions from $\mathbb{R}^n$ to $\mathbb{R}^n$ that is parameterized by some index set $\Lambda$ and $\sigma : [0, \infty) \to \Lambda$ is the switching signal which is a piecewise constant function of time. At each instant of time, the switching signal may be an arbitrary function of time or system states $x(t)$, or even might be user specified. It is assumed that the state of (7) does not jump at the switching instants, i.e. the solution of $x(t)$ is continuous everywhere. Also the case of infinitely fast switching (chattering) is not considered in the analysis. In the case where the individual subsystems of the switched system are linear, equation (7) reduces to (8).

$$\dot{x} = \Psi_\sigma x \quad (8)$$

The stability of switched systems, where the switching laws are not known a-priori, is not a trivial problem. A necessary condition for stability under arbitrary switching is that the individual subsystems must be stable [20]. However for certain types of switching signals it is possible to get instability even though all the subsystems are asymptotically stable [20].

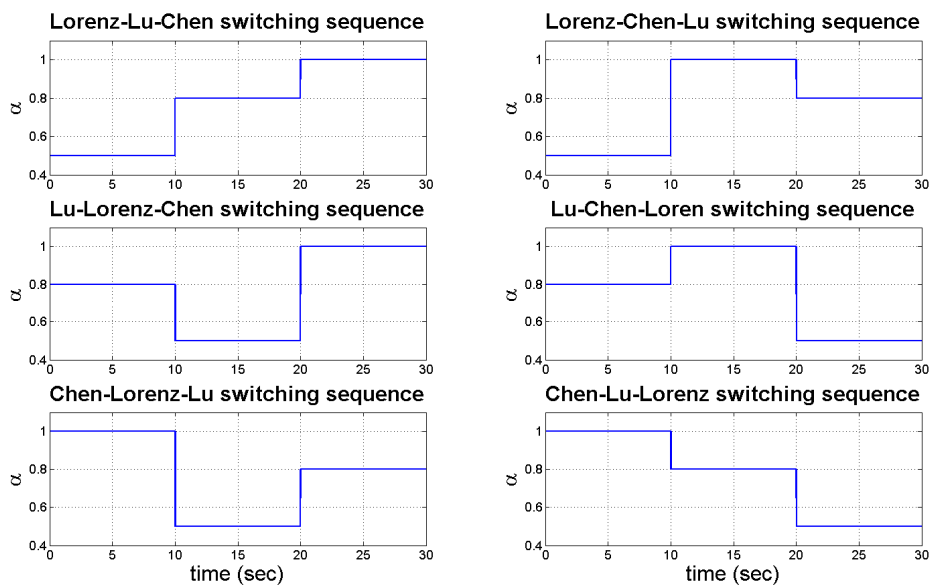

**Figure 1: Switching law for parameter $\alpha$ in the unified chaotic system.**

If sudden switching in the parameter $\alpha$ is considered in the unified chaotic system, the resulting system can be viewed as a switched chaotic system which is a new concept, studied in this paper. We have also shown the synchronization in two such switched chaotic systems using a stabilizing nonlinear state feedback control law. Before studying synchronization, the dynamical behaviors of switched chaotic systems are explored in a more detailed manner. Figure 1 shows the switching instants of the parameter $\alpha$, so that the unified chaotic system jumps between Lorenz-Chen-Lu family. In fact, these three subsystems (or family of attractors) can be switched in different combinations to produce different switched chaotic systems. In the contemporary literatures, synchronization and control performance is studied amongst various class of chaotic and hyper-chaotic systems [22], [23] along with recent developments of optimization based generalized control scheme design for special class of chaotic systems e.g. the multi-wing Lorenz family [39]. But to the best of the authors' knowledge there is no work aiming to the design of a global controller to ensure guaranteed synchronization amongst different switching conditions of the unified chaotic system family. In Figure 1 the piecewise continuous switching signal $\sigma(t)$ changes the value of $\alpha$ in the unified chaotic system model of (1) thereby switching between different chaotic models. Thus $\sigma(t) \in [0,1]$ which is the range where $\alpha$ can vary in the Unified chaotic system (1). During the first 10 seconds the system behaves like either of the Lorenz/Lu/Chen system in Figure 1. In the next interval (10-20 sec) the system evolves as a different chaotic system other than the starting one and in the final phase (20-30 sec) the chaotic system again evolves as a different one. Such a configuration makes the whole switched chaotic system as a combination of three chaotic subsystems which is more complex than the usual ones, where synchronization is generally targeted within two fixed structure (similar or different) chaotic systems. In this paper, the master and slave systems are assumed to be driven by the same switching law.

## 3. Synchronization of the switched chaotic systems by a nonlinear state feedback controller

### 3.1 *Chaos synchronization as a stabilization problem*

For the nonlinear error dynamical system given by (6) the feedback control law can be designed to have a linear and a nonlinear component, i.e.

$$
\begin{aligned}
\dot{e} &= Ae + f(u) \\
&= Ae + f_{nonlinear}(u) + f_{linear}(u) \\
&= Ae + U_1 e + U_2 e
\end{aligned}
\tag{9}
$$

where, $U_1$ and $U_2$ are matrices of appropriate dimensions, which contribute to the nonlinear and linear parts of the control input respectively.

The feedback control law is chosen as in (10).

$$
u = \begin{bmatrix} 0 \\ z_m e_1 + x_s e_3 \\ -y_m e_1 - x_s e_2 \end{bmatrix} + \begin{bmatrix} k_1 e_1 \\ k_2 e_2 \\ k_3 e_3 \end{bmatrix}
\tag{10}
$$

i.e.

$$
u = \begin{bmatrix} 0 & 0 & 0 \\ z_m & 0 & x_s \\ -y_m & -x_s & 0 \end{bmatrix} \begin{bmatrix} e_1 \\ e_2 \\ e_3 \end{bmatrix} + \begin{bmatrix} k_1 & 0 & 0 \\ 0 & k_2 & 0 \\ 0 & 0 & k_3 \end{bmatrix} \begin{bmatrix} e_1 \\ e_2 \\ e_3 \end{bmatrix}
\tag{11}
$$

Using (9) and (11) we get

$$
\dot{e} = (A + U_1) e + U_2 e
\tag{12}
$$

$$
\Rightarrow \dot{e} = \tilde{A} e + BK e
\tag{13}
$$

where, $\tilde{A} = A + U_1 = \begin{bmatrix} -(25\alpha + 10) & (25\alpha + 10) & 0 \\ (28 - 35\alpha) & (29\alpha - 1) & 0 \\ 0 & 0 & -\left(\dfrac{\alpha + 8}{3}\right) \end{bmatrix}$, $U_2 = BK$, $B = \begin{bmatrix} 1 \\ 1 \\ 1 \end{bmatrix}$,

$K = \begin{bmatrix} k_1 & k_2 & k_3 \end{bmatrix}$.

Here, $K$ represents the state feedback controller gains and $B$ is essentially the distribution matrix which gives a weighting on these gains and adds these control signals to the individual states of the error dynamical system.

**Remark 1:**

For a given set of controller gains $K$, the error dynamical system (13) can be treated as an uncertain system with the variable parameter $\alpha$. Here the values of $\tilde{A}$ are uncertain but belongs to a polytopic set $\Omega = Cov\left[\tilde{A}_1, \tilde{A}_2, \tilde{A}_3 \ldots \ldots, \tilde{A}_{N_m}\right]$;

where, $Cov$ refers to a convex hull.

Thus any $\tilde{A} \in \Omega$ can be expressed as

$$\hat{A} = \sum_{i=1}^{N_m} w_i(\alpha)[\tilde{A}_i] \qquad (14)$$

Here, $N_m$ is the number of multiple models, and $w_i$ the weighting functions constrained between 0 and 1 which satisfies

$$\sum_{i=1}^{N_m} w_i(\alpha) = 1, \ \ \forall \alpha \in [0,1] \qquad (15)$$

Generally, we have to search for a solution to the matrix inequalities, all over the uncertain system $\Omega$. However, due to the properties of the polytopic system, solution can be found only at the polytopic vertices instead of all points within the polytope. The polytope vertices are found by taking all possible combinations of the maximum and minimum values of each element of $\tilde{A}$ containing $\alpha$. Since there are 5 such elements, hence in this problem we have 32 different cases of $\tilde{A}$, i.e., $\Omega = Cov\left[\tilde{A}_1, \tilde{A}_2, \tilde{A}_3 \ldots \ldots, \tilde{A}_{32}\right]$.

**Theorem 2:** [40]
The following two statements in (16) and (17) are equivalent:
 (a) The switched linear system

$$\dot{x}(t) = A_{\sigma(t)} x(t) \qquad (16)$$

 where, $A_{\sigma(t)} \in \{A_1, A_2, \ldots, A_N\}$, is asymptotically stable under arbitrary switching.

 (b) The linear time varying system

$$\dot{x}(t) = A(t) x(t) \qquad (17)$$

 where, $A(t) \in \text{Cov}\{A_1, A_2, \ldots, A_N\}$, is robustly asymptotically stable.

The proof is excluded here for brevity. Please look in [40] and the references therein.

**Corollary 2.1:**
The following two statements in (18) and (19) are equivalent:
 (a) The switched linear system $A_{\sigma(t)}$ with the state feedback control $BK$

$$\dot{x}(t) = \left(A_{\sigma(t)} + BK\right) x(t) \qquad (18)$$

 where, $A_{\sigma(t)} \in \{A_1, A_2, \ldots, A_N\}$, is asymptotically stable under arbitrary switching.

 (b) The linear time varying system $A(t)$ with the state feedback control $BK$

$$\dot{x}(t) = \left(A(t) + BK\right) x(t) \qquad (19)$$

 where, $A(t) \in \text{Cov}\{A_1, A_2, \ldots, A_N\}$, is robustly asymptotically stable.

It is important to notice that for both the switched system (18) and linear time varying system (19), one part containing the term $BK$ is constant whereas the other part changes over time. Equation (18) investigates the stability for switching between the vertices of the switching combinations $A_{\sigma(t)}$. Whereas equation (19) says if the system $A(t)$ is stable at the vertices of the convex hull then the stability is automatically satisfied within the internal points of the convex region. The proof of this statement and the corresponding controller design can be considered as the main contributions of the paper.

**Proof that statement 2.1(a) implies 2.1(b):**

For stability of the system in (18), there must exist a common Lyapunov function which decreases along the trajectory of the system, as the system switches arbitrarily between the individual subsystems.

For proving that the statement in Corollary 2.1(a) implies 2.1(b) we have the common quadratic Lyapunov function $V(x)$ of the form

$$V(x) := x^T P x \tag{20}$$

where, $P$ is a symmetric positive definite matrix of appropriate dimensions, i.e.

$$P > 0 \tag{21}$$

we must have $\dot{V}(x) < 0$ to ensure stability.

Now,

$$\begin{aligned}
\dot{V}(x) &= \dot{x}^T P x + x^T P \dot{x} \\
&= \dot{x}^T \left\{ \left( A_{\sigma(t)} + BK \right)^T P + P \left( A_{\sigma(t)} + BK \right) \right\} x
\end{aligned} \tag{22}$$

Hence for stability of the system (18) we must have

$$\left( A_{\sigma(t)} + BK \right)^T P + P \left( A_{\sigma(t)} + BK \right) < 0 \tag{23}$$

i.e.

$$\left( A_i + BK \right)^T P + P \left( A_i + BK \right) < 0 \quad \forall \ i \in \{1, \cdots, N\}. \tag{24}$$

Now for the proof, it is necessary to show that the robust stability of system (19) is also satisfied by the same conditions (21) and (24).

Since $A(t)$ in (19) belongs to $\mathrm{Cov}\{A_1, A_2, \ldots, A_N\}$, it is possible to express $A(t)$ as a convex combination of the vertices, i.e.

$$A(t) = \sum_{i=1}^{N} w_i A_i, \quad \text{where,} \sum_{i=1}^{N} w_i = 1 \tag{25}$$

Now the same Lyapunov function (20) should decrease monotonically for robust stability of (19). Thus condition (26) must hold.

$$(A(t) + BK)^T P + P(A(t) + BK) < 0 \tag{26}$$

Using (25), the left hand side of the inequality (26) can be written as:

$$\left(\sum_{i=1}^{N} w_i A_i + BK\right)^T P + P\left(\sum_{i=1}^{N} w_i A_i + BK\right)$$

$$= w_1 A_1^T P + w_2 A_2^T P + \ldots + w_N A_N^T P + P w_1 A_1 + P w_2 A_2 + \ldots P w_N A_N + (BK)^T P + P(BK)$$

$$= w_1\left(A_1^T P + PA_1\right) + w_2\left(A_2^T P + PA_2\right) + \ldots + w_N\left(A_N^T P + PA_N\right) + (BK)^T P + P(BK)$$

$$= w_1\left\{(A_1 + BK)^T P + P(A_1 + BK)\right\} + w_2\left\{(A_2 + BK)^T P + P(A_2 + BK)\right\} + \ldots$$

$$\qquad + w_N\left\{(A_N + BK)^T P + P(A_N + BK)\right\} - \left\{w_1(BK)^T P + w_1 P(BK)\right\}$$

$$\qquad - \left\{w_2(BK)^T P + w_2 P(BK)\right\} - \ldots - \left\{w_N(BK)^T P + w_N P(BK)\right\} + (BK)^T P + P(BK)$$

$$= w_1\left\{(A_1 + BK)^T P + P(A_1 + BK)\right\} + w_2\left\{(A_2 + BK)^T P + P(A_2 + BK)\right\} + \ldots$$

$$\qquad + w_N\left\{(A_N + BK)^T P + P(A_N + BK)\right\} - (BK)^T P\left[w_1 + w_2 + \ldots + w_{N_m} - 1\right]$$

$$\qquad - PBK\left[w_1 + w_2 + \ldots + w_{N_m} - 1\right]$$

$$= w_1\left\{(A_1 + BK)^T P + P(A_1 + BK)\right\} + w_2\left\{(A_2 + BK)^T P + P(A_2 + BK)\right\} + \ldots$$

$$\qquad + w_N\left\{(A_N + BK)^T P + P(A_N + BK)\right\} \qquad \left[\because \sum_{i=1}^{N} w_i = 1\right] \tag{27}$$

$$= \sum_{i=1}^{N} w_i\left\{(A_i + BK)^T P + P(A_i + BK)\right\}$$

Since the condition (24) holds, therefore in (27), all the terms of the form $\left\{(A_i + BK)^T P + P(A_i + BK)\right\} < 0, \forall i \in \{1, \cdots, N\}$. Also each of these terms are multiplied by $w_i \in [0,1]$, hence,

$$w_1\left\{(A_1 + BK)^T P + P(A_1 + BK)\right\} + w_2\left\{(A_2 + BK)^T P + P(A_2 + BK)\right\} + \ldots$$

$$\qquad + w_N\left\{(A_N + BK)^T P + P(A_N + BK)\right\} < 0 \tag{28}$$

$$\Rightarrow (A(t) + BK)^T P + P(A(t) + BK) < 0$$

which was required to hold in (26) for the proof. $\square$

**Proof that statement 2.1(b) implies 2.1.(a):**

Now in order to prove the reverse way i.e. for proving that statement in Corollary 2.1(b) implies 2.1(a) we have the system in (19). For robust stability of (19), let the common quadratic Lyapunov function have the same form as in (20) where $P > 0$. To ensure stability of (19) we must have $\dot{V}(x) < 0$.

Now,

$$\dot{V}(x) = \dot{x}^T P x + x^T P \dot{x}$$

$$= \left[ \left( A(t) + BK \right) x(t) \right]^T P x(t) + x^T(t) P \left( A(t) + BK \right) x(t) \qquad (29)$$

$$= x^T(t) \left\{ \left( A(t) + BK \right)^T P + P \left( A(t) + BK \right) \right\} x(t)$$

Hence for stability of the system (19), equation (30) must hold.

$$\left( A(t) + BK \right)^T P + P \left( A(t) + BK \right) < 0 \qquad (30)$$

In order to prove the equivalence of the two statements in Corollary 2.1(b) to 2.1(a), it is necessary that the same $P$ matrix should also satisfy the stability of the switched system in (18). Now, using (25) from expression (30) we have:

$$\left( \sum_{i=1}^{N} w_i A_i + BK \right)^T P + P \left( \sum_{i=1}^{N} w_i A_i + BK \right) < 0$$

$$\Rightarrow \sum_{i=1}^{N} w_i A_i^T P + \left( BK \right)^T P + P \sum_{i=1}^{N} w_i A_i + P \left( BK \right) < 0$$

$$\Rightarrow \sum_{i=1}^{N} w_i \left( A_i^T P + PA \right) + \left( BK \right)^T P + P \left( BK \right) < 0$$

$$\Rightarrow \sum_{i=1}^{N} w_i \left\{ \left( A_i + BK \right)^T P + P \left( A_i + BK \right) \right\} - \sum_{i=1}^{N} w_i \left\{ \left( BK \right)^T P + P \left( BK \right) \right\} + \left( BK \right)^T P + P \left( BK \right) < 0$$

$$\Rightarrow \sum_{i=1}^{N} w_i \left\{ \left( A_i + BK \right)^T P + P \left( A_i + BK \right) \right\} + \left\{ \left( BK \right)^T P + P \left( BK \right) \right\} \left( 1 - \sum_{i=1}^{N} w_i \right) < 0$$

$$\Rightarrow \sum_{i=1}^{N} w_i \left\{ \left( A_i + BK \right)^T P + P \left( A_i + BK \right) \right\} < 0 \quad \left[ \because \sum_{i=1}^{N} w_i = 1 \right]$$

$$(31)$$

Although (31) is a summation term for all the cases ($i = 1, \cdots, N$) and it is also less than zero where the weights $w_i \in [0,1]$, it is not straightforward to conclude that the individual terms $\left\{ \left( A_i + BK \right)^T P + P \left( A_i + BK \right) \right\}, \forall i$ are also less than zero. The converse problem is generally easy to prove as done in the forward proof i.e. if $\sum_{i=1}^{N} w_i = 1, w_i > 0$ and $\left\{ \left( A_i + BK \right)^T P + P \left( A_i + BK \right) \right\} < 0, \forall i$, therefore

$$\sum_{i=1}^{N} w_i \left\{ \left( A_i + BK \right)^T P + P \left( A_i + BK \right) \right\} < 0. \qquad (32)$$

This is because each term within the curly bracket in (32) is multiplied by a positive fraction and the overall expression will also be less than zero. For the reverse proof, however, we note that since, equation (31) is valid for all values of $w_i \in [0,1]$, we get:

$$\left(A_1 + BK\right)^T P + P\left(A_1 + BK\right) < 0 \quad \text{by setting } w_1 = 1, w_i = 0, \forall\, i \neq 1$$

$$\left(A_2 + BK\right)^T P + P\left(A_2 + BK\right) < 0 \quad \text{by setting } w_2 = 1, w_i = 0, \forall\, i \neq 2$$

$$\vdots$$

$$\left(A_N + BK\right)^T P + P\left(A_N + BK\right) < 0 \quad \text{by setting } w_N = 1, w_i = 0, \forall\, i \neq N$$

(33)

The expressions in (33) can be written concisely as (34).

$$\left(A_i + BK\right)^T P + P\left(A_i + BK\right) < 0, \forall\, i \in \left\{1, \cdots, N\right\}$$

(34)

This is the stability criteria for the switched system $\left(A_{\sigma(t)} + BK\right)$ with a common Lypunov function $P$ as obtained from equations (23) and (24) which was required for the proof. □

Therefore, we here proved from both ways that the expressions in Corollary 2.1(a) and (b) are equivalent. This implies that stability of a switched linear system with the state feedback controller, under arbitrary switching at the polytope vertices, and the robust stability of the corresponding linear time varying system with the same state feedback controller are indeed equivalent.

**Theorem 3:**

For $\tilde{A}$ with polytopic type uncertainty, stability of the system (13) is satisfied if $\left(\tilde{A}_i + BK\right)^T P + P(\tilde{A}_i + BK) < 0$, where $P > 0$ is a positive definite matrix of appropriate dimensions.

**Proof:**

Let $V\left(e\right)$ be a common quadratic Lyapunov function of the form

$$V\left(e\right) := e^T P e$$

(35)

where, $P$ is a symmetric positive definite matrix of appropriate dimensions, i.e.

$$P > 0$$

(36)

The closed loop system with the feedback controller as in (9), becomes an autonomous system. From (13) the autonomous error dynamical system can be written as

$$\dot{e} = \left(\tilde{A} + BK\right)e$$

(37)

Considering that $\alpha$ is a variable in $\tilde{A}$, which changes according to the switching signal $\sigma$, it is possible to write $\tilde{A}$ in the generalized form as $\hat{A}$ as shown in (14).

Thus (37) becomes

$$\dot{e} = \left(\hat{A} + BK\right)e$$

(38)

For Lyapunov stability of the system (38), $\dfrac{\partial V\left(e\right)}{\partial t} < 0$ must hold so that $V\left(e\right)$ is monotonically decreasing along the trajectory of the system.

Thus,

$$\frac{\partial V(e)}{\partial t} = \dot{e}^T Pe + e^T P\dot{e}$$

$$= e^T \left(\hat{A} + BK\right)^T Pe + e^T P\left(\hat{A} + BK\right)e \qquad (39)$$

$$= e^T \left\{\left(\hat{A} + BK\right)^T P + P\left(\hat{A} + BK\right)\right\}e$$

Thus for a monotonically decreasing Lyapunov function (40) must hold.

$$\left(\hat{A} + BK\right)^T P + P\left(\hat{A} + BK\right) < 0 \qquad (40)$$

Hence the LMI in (40) must be satisfied for the whole region $\Omega$. However, since it is a polytopic type uncertainty, only the vertices of the convex hull of $\Omega$ need to be satisfied for stability. Thus (40) reduces to (41) under this condition.

$$\left(\tilde{A}_i + BK\right)^T P + P(\tilde{A}_i + BK) < 0 \qquad (41)$$

where $i \in \{1, \cdots, N_m\}$

Thus for (36) we have $V(e) > 0$ and for (41) we have $\dot{V}(e) < 0$. Thus $\|e(t)\| \rightarrow 0$ and hence the system is stable. □

However considering (41) and (18), (19) we find that the system stability is ensured only if it arbitrarily switches at the polytope vertices. Since $\alpha$ of the unified chaotic system in (1) can take any value in the interval $[0,1]$ it is necessary to show that (41) is a sufficient condition for stability, for all convex combinations of the polytope vertices $\tilde{A}_i$. In other words switching at the vertices as well as inside the polytope is not covered by (18) and (19) which has been proved next. Although it is trivial in the stabilization problem of switched linear systems [14], but for the sake of completeness and also to establish the argument in the case of nonlinear chaotic system synchronization with arbitrary switching, we reiterate the mathematical formulation and its proof in the present context.

**Theorem 4** [14]:

Stability of the switched system $\tilde{A}_i \ \forall i \in \{1, \cdots, N_m\}$ with the state feedback controller $BK$, as given by (41) is a sufficient condition for the stability of arbitrary switching between all convex combinations of $\tilde{A}_i$.

**Proof:**

From (14) and (15) any convex combination of $\tilde{A}_i$ can be expressed as

$$\hat{A} = \sum_{i=1}^{N_m} w_i \tilde{A}_i \qquad (42)$$

where $\sum_{i=1}^{N_m} w_i = 1$.

Hence if $\hat{A}$ with the feedback controller, has to be stable, given condition (41), it must satisfy the Lyapunov stability condition with the same $P$ matrix as in (41). Thus the following condition in (43) must hold.

$$\left(\hat{A}+BK\right)^T P + P\left(\hat{A}+BK\right) < 0 \qquad (43)$$

Using (42) the left hand side of inequality (43) can be written as,

$$\left(\sum_{i=1}^{N_m} w_i \tilde{A}_i + BK\right)^T P + P\left(\sum_{i=1}^{N_m} w_i \tilde{A}_i + BK\right)$$

$$\Rightarrow \left(\sum_{i=1}^{N_m} w_i \tilde{A}_i\right)^T P + \left(BK\right)^T P + P\left(\sum_{i=1}^{N_m} w_i \tilde{A}_i\right) + P\left(BK\right)$$

$$\Rightarrow \left(w_1 \tilde{A}_1^T + w_2 \tilde{A}_2^T + \cdots + w_{N_m} \tilde{A}_{N_m}^T\right) P + P\left(w_1 \tilde{A}_1 + w_2 \tilde{A}_2 + \cdots + w_{N_m} \tilde{A}_{N_m}\right) + \left(BK\right)^T P + P\left(BK\right)$$

$$\Rightarrow w_1\left\{\left(\tilde{A}_1 + BK\right)^T P + P\left(\tilde{A}_1 + BK\right)\right\} + w_2\left\{\left(\tilde{A}_2 + BK\right)^T P + P\left(\tilde{A}_2 + BK\right)\right\} + \cdots$$

$$\quad + w_{N_m}\left\{\left(\tilde{A}_{N_m} + BK\right)^T P + P\left(\tilde{A}_{N_m} + BK\right)\right\}$$

$$\quad + \left(BK\right)^T P\left\{1 - w_1 - w_2 \cdots - w_{N_m}\right\} + P\left(BK\right)\left\{1 - w_1 - w_2 \cdots - w_{N_m}\right\}$$

$$\Rightarrow w_1\left\{\left(\tilde{A}_1 + BK\right)^T P + P\left(\tilde{A}_1 + BK\right)\right\} + w_2\left\{\left(\tilde{A}_2 + BK\right)^T P + P\left(\tilde{A}_2 + BK\right)\right\} + \cdots$$

$$\quad + w_{N_m}\left\{\left(\tilde{A}_{N_m} + BK\right)^T P + P\left(\tilde{A}_{N_m} + BK\right)\right\} \quad \left[\because 1 - \sum_{i=1}^{N_m} w_i = 0\right]$$

(44)

From (41) we find that each individual polytope vertex $\tilde{A}_i \ \forall i \in \{1, \cdots, N_m\}$ is stable with feedback gain $BK$. Thus in (44) all the terms of the form $\left(\tilde{A}_i + BK\right)^T P + P(\tilde{A}_i + BK) \ \forall i \in \{1, \cdots, N_m\}$ is less than zero. Also each of these terms are multiplied by $w_i \in [0,1]$, hence,

$$w_1\left\{\left(\tilde{A}_1 + BK\right)^T P + P\left(\tilde{A}_1 + BK\right)\right\} + w_2\left\{\left(\tilde{A}_2 + BK\right)^T P + P\left(\tilde{A}_2 + BK\right)\right\} + \cdots$$

$$\quad + w_{N_m}\left\{\left(\tilde{A}_{N_m} + BK\right)^T P + P\left(\tilde{A}_{N_m} + BK\right)\right\} < 0 \qquad (45)$$

$$\Rightarrow \left(\hat{A}+BK\right)^T P + P\left(\hat{A}+BK\right) < 0$$

Which was required to hold in (43) for the proof. □

## 3.2 BMI formulation for synchronizing controller design and reduction in LMIs

The controller design problem is to find the values of the feedback gains $K$ such that the matrix inequalities (41) hold for some values of $P > 0$. However since (41) contains both unknown variables $K$ and $P$ multiplied together, hence they are in the form of bilinear matrix inequalities or BMIs and cannot be solved using conventional

LMI framework [30]. Thus the inequalities need to be recast suitably to facilitate the numerical solution by standard LMI solvers.

**Theorem 5:**
The controller design problem satisfying (41) is solvable using standard LMI framework by recasting it into the following form [37]:

$$\tilde{A}_i Y + Y \tilde{A}_i^T + B_i K_a + K_a^T B_i^T < 0 \qquad (46)$$

where, $K_a := KY$ and $Y > 0 , Y := P^{-1}$

**Proof:**
Since

$$\begin{aligned}
& P > 0 \\
& \Rightarrow P^{-1} P P^{-1} > 0 \\
& \Rightarrow P^{-1} > 0 \left[ P^{-1} P = I \right] \\
& \Rightarrow Y > 0,
\end{aligned} \qquad (47)$$

therefore from (41) we have

$$\begin{aligned}
& \left( \tilde{A}_i + BK \right)^T P + P(\tilde{A}_i + BK) < 0 \\
& \Rightarrow (\tilde{A}_i^T + K^T B^T) P + P(\tilde{A}_i + BK) < 0 \\
& \Rightarrow (\tilde{A}_i^T + K^T B^T) Y^{-1} + Y^{-1}(\tilde{A}_i + BK) < 0 \quad [\because Y^{-1} := P] \\
& \Rightarrow \tilde{A}_i^T Y^{-1} + K^T B^T Y^{-1} + Y^{-1} \tilde{A}_i + Y^{-1} BK < 0 \\
& \Rightarrow Y \tilde{A}_i^T + YK^T B^T + \tilde{A}_i Y + BKY < 0 \quad [\text{pre \& post-multiplying by } Y] \\
& \Rightarrow Y \tilde{A}_i^T + \tilde{A}_i Y + Y(K_a Y^{-1})^T B^T + BK_a Y^{-1} Y < 0, \quad [\text{where, } K_a := KY] \\
& \Rightarrow Y \tilde{A}_i^T + \tilde{A}_i Y + Y(Y^{-1})^T K_a^T B^T + BK_a < 0 \\
& \Rightarrow Y \tilde{A}_i^T + \tilde{A}_i Y + Y(Y^{-1}) K_a^T B^T + BK_a < 0 \quad [\text{since, } (Y^{-1})^T = Y^{-1} (\text{symmetric})] \\
& \Rightarrow Y \tilde{A}_i^T + \tilde{A}_i Y + K_a^T B^T + BK_a < 0. \quad \square
\end{aligned} \qquad (48)$$

Thus (47) and (48) are in the form of LMIs with variables $Y$ and $K_a$.
Finally the state feedback controller $K$ can be obtained by using the following relationship.

$$K = K_a Y^{-1} \qquad (49)$$

# 4. Numerical simulation examples

For the controller design procedure as illustrated in Section 3.2, (47) and (48) are solved as a feasibility problem using YALMIP Toolbox [13]. From the LMI solving results, the feedback gain is obtained as $K_a = \begin{bmatrix} -0.3012 & -0.6857 & 0.5681 \end{bmatrix}$ and the values of $Y$ matrix and $P$ matrix are as follows:

$$Y = \begin{bmatrix} 0.0385 & -0.0014 & 0.0112 \\ -0.0014 & 0.0013 & -0.0005 \\ 0.0112 & -0.0005 & 0.5753 \end{bmatrix} \text{ and } P = \begin{bmatrix} 27.2002 & 28.9674 & -0.5062 \\ 28.9674 & 779.2615 & 0.0738 \\ -0.5062 & 0.0738 & 1.7481 \end{bmatrix}.$$

The eigenvalues of the obtained $P$ matrix are as follows:

$$\text{eig}(P) = \begin{bmatrix} 1.7375 & 26.0968 & 780.3756 \end{bmatrix}.$$

Since all the eigenvalues are positive, hence $P$ is positive definite and the inequality conditions in (41) and (36) are satisfied.

Next the actual feedback gain $K = \begin{bmatrix} k_1 & k_2 & k_3 \end{bmatrix}$ of the nonlinear controller in (11) is obtained as

$$K = K_a Y^{-1} = K_a P = \begin{bmatrix} -28.3433 & -543.0217 & 1.0950 \end{bmatrix} \tag{50}$$

Using the controller gains obtained from LMI solvers i.e. (50) we now show the synchronization performance amongst different configurations of switching in the unified chaotic system family in Figure 2. Here, the corresponding switching laws i.e. switching sequences of parameter $\alpha$ are shown in Figure 1. Depending on the starting nature and intermediate type of switching amongst these three families of chaotic systems, different dynamics of the states are observed as shown in Figure 2. It is noticeable that the single controller is capable of faithfully enforcing the slave trajectory to follow the master trajectory for all the six cases in Figure 2 where the switching sequences of the key parameter $\alpha$ of the unified chaotic system family are shown in Figure 1.

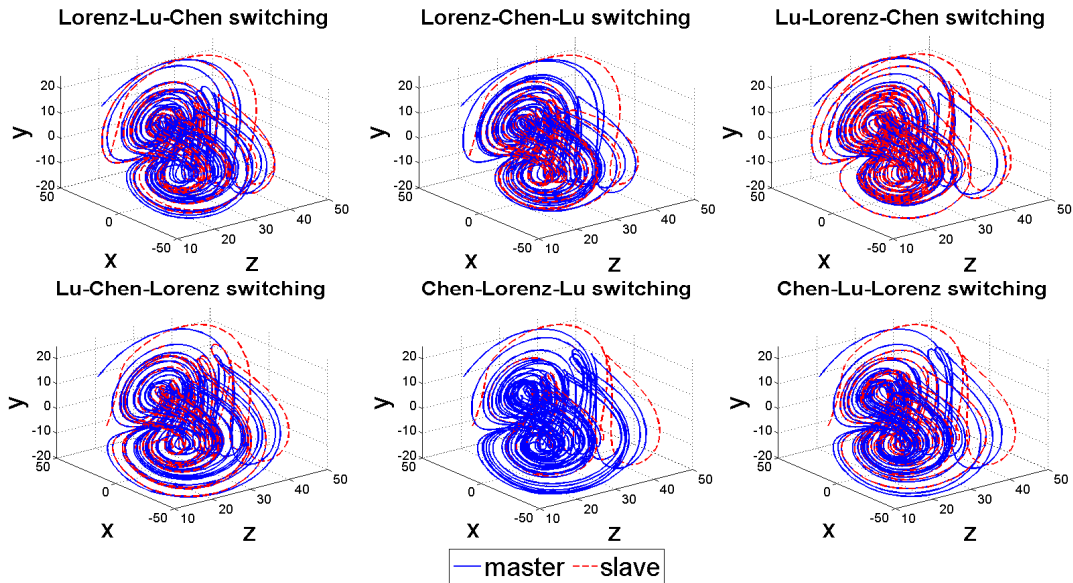

**Figure 2: Synchronization in different configurations of switched unified chaotic systems.**

For all the simulations reported here, the initial conditions of the master and slave chaotic systems are chosen as $(x_{m0}, y_{m0}, z_{m0}) = (15, 20, 10)$ and $(x_{s0}, y_{s0}, z_{s0}) = (25, -5, 15)$ respectively. The synchronization error is defined as the $L_2$-norm of the error vector as shown in [41], i.e.

$$e = \sqrt{e_x^2 + e_y^2 + e_z^2} \qquad (51)$$

For efficient synchronization the error should decay with time and its final value should tend towards a small finite value. With the same controller gains designed via satisfying the LMI conditions in (50) we can achieve guaranteed synchronization amongst various switching conditions of the unified chaotic system which are explored next.

### 4.1 Sudden step switching in key parameter α

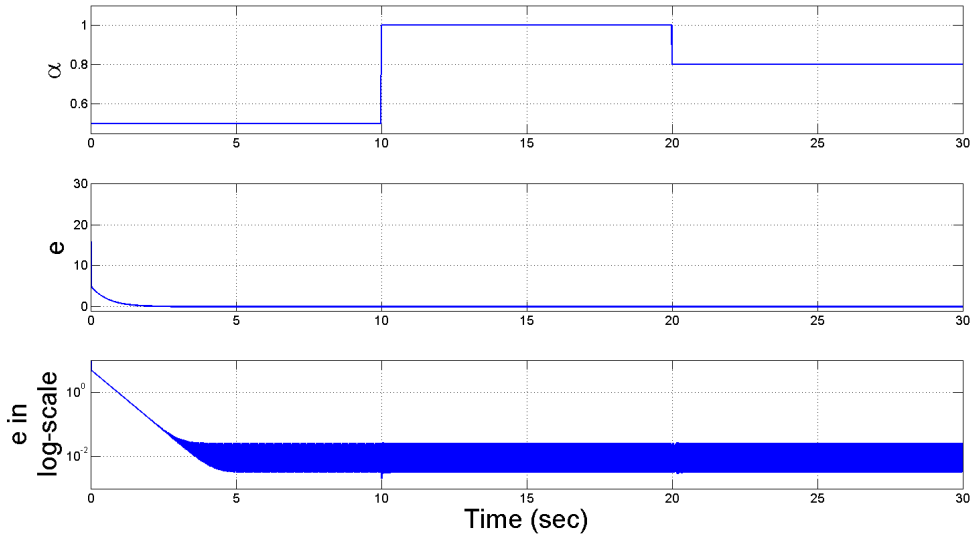

**Figure 3: Sudden switching law for α and the corresponding synchronization error in normal and logarithmic scale.**

Here we consider that the key parameter $\alpha$ of the error dynamical system (6) and hence the type of master and slave chaotic systems are suddenly switched as a step function. It is interesting to note from the time evolution of the synchronization error in Figure 3 and the time evolution of the states of the master-slave systems in Figure 4 that even though the switching occurs (i.e. the value of $\alpha$ changes) at the $t = 10$ and $t = 20$ sec, the slave states are perfectly synchronized with that of the master, due to the effectiveness of the designed global nonlinear controller. However, observation of the patterns of the state evolution during the three time intervals (i.e. 0-10 sec, 10-20 sec and 20-30 sec) in Figure 4 shows clearly that they represent different chaotic systems.

Also the time required for synchronization of the master and slave states is very less. For the $z$ state, there is a perceptible time for the synchronization process as is evident from Figure 4. For the $x$ and $y$ states, the synchronization is almost instantaneous. Also, in order to show the goodness of synchronization of the states of master and slave, the time evolution of the error is shown in normal and logarithmic scale in Figure 3. Figure 5 shows the master-slave synchronization in the 2D and 3D phase space diagrams using the proposed global controller when the initial conditions of the two unified chaotic systems are different. Due to the nonlinear terms in the controller while satisfying the Lyapunov stability criterion for switched systems, the other possible switching combinations among the unified chaotic system family apart from that shown in Figure 3 and Figure 5 can also be synchronized efficiently as explored in Figure 1. The corresponding master-slave synchronization performances and the evolution of the state trajectories with the controller have already been shown in Figure 2.

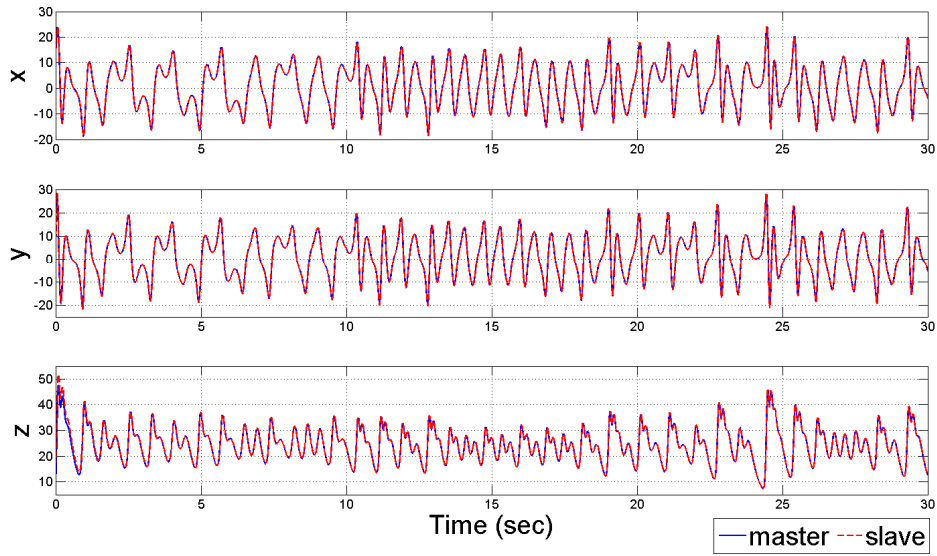

Figure 4: Time evolution of the master-slave states for sudden switching in *α*.

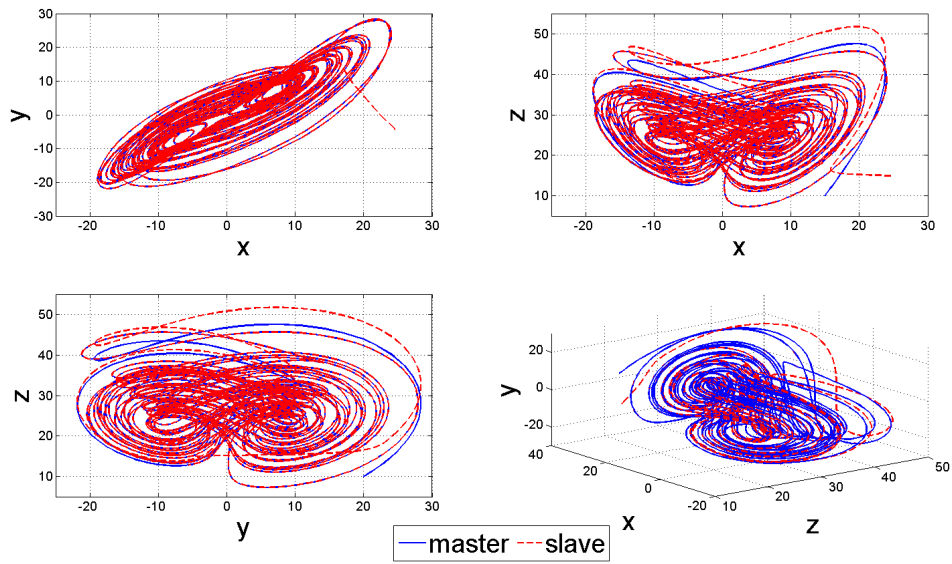

Figure 5: Phase space representation of the master-slave states for sudden step change switching in *α*.

## 4.2 Sinusoidal switching in key parameter α

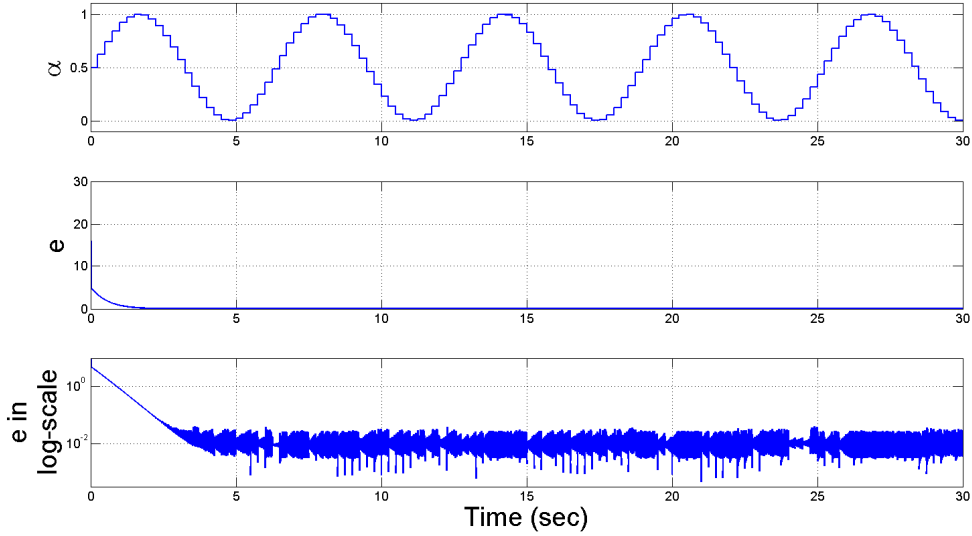

**Figure 6: Sinusoidal switching law for α and the corresponding synchronization error in normal and logarithmic scale.**

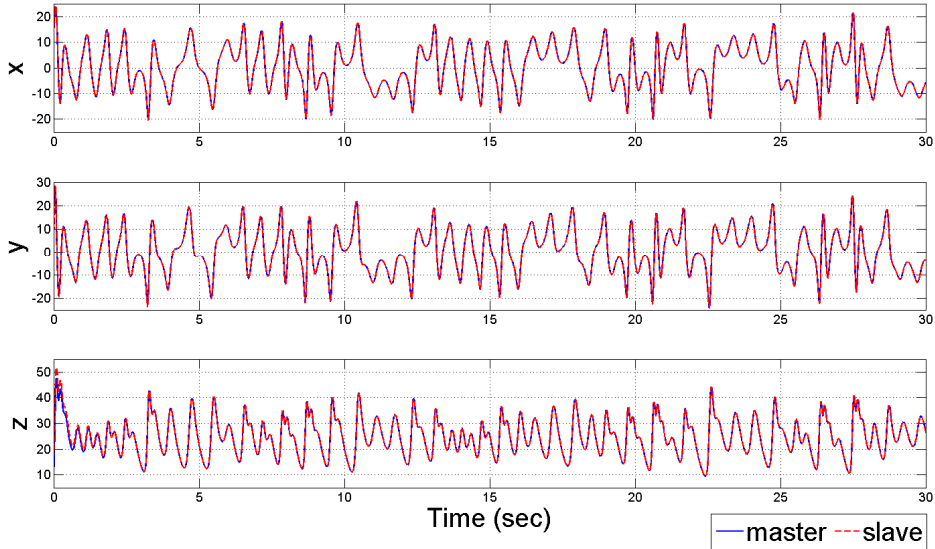

**Figure 7: Time evolution of the master-slave states for sinusoidal switching in α.**

Next, the switching law for the key parameter (α) of the unified chaotic system is defined as a sinusoidal signal along with a sample and zero-order-hold (ZOH) operation. The sampling time for the switching law is chosen as 0.25 sec. The sinusoidal switching law is having a frequency of 1 rad/sec with an amplitude of 0.5 and a bias of 0.5, so that $\alpha \in [0,1]$ in all cases representing different dynamics among the unified chaotic system family. The key parameter $\alpha$ here continuously switches to take different values within the unified chaotic system family and the synchronization performance is shown in Figure 6. It is noticeable from Figure 6 that there exist small fluctuations of the order of $10^{-2}$ in log-scale, due the presence of more number of switching in the chaotic systems in comparison with the case in Figure 1 and Figure 3.

Even though the master and slave systems are continuously getting switched multiple times, it is interesting to note that the states once synchronized do not diverge again to a great extent due to the strong nonlinear structure and the effective LMI based design technique of the synchronizing controller as shown in Figure 7. The corresponding synchronization performance in the phase space diagram has been in Figure 8.

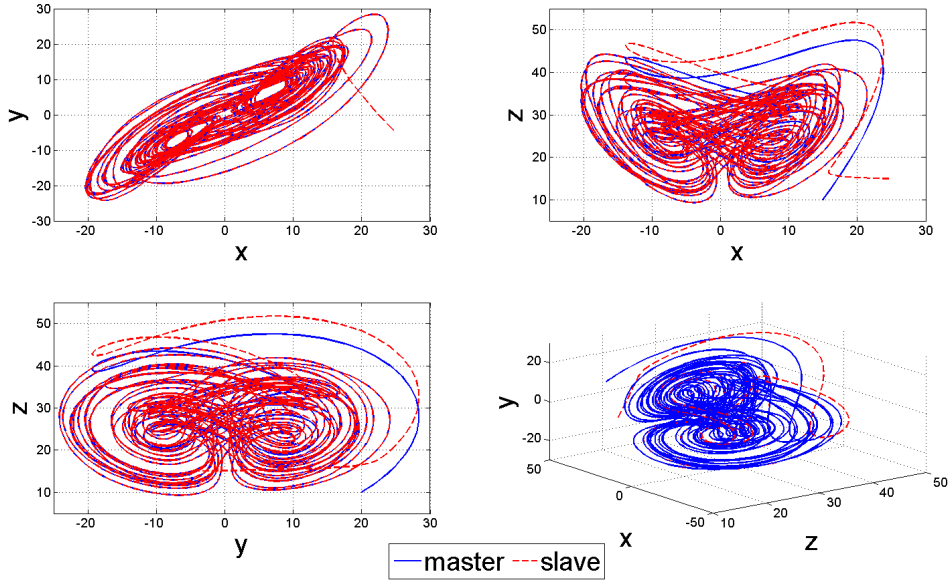

**Figure 8: Phase space representation of the master-slave states for sinusoidal switching in $\alpha$.**

### 4.3  Chirp signal switching in key parameter α

The switching law of the key parameter of the unified chaotic system is then varied as a chirp signal, along with a chosen sampling operation after every 0.1 sec and a ZOH operation. The frequency of the chirp signal has been increased from 0.1 to 1 Hz within 30 sec of simulation time. The generated chirp signal is halved, passed through the sample and hold operations and then added with a bias of 0.5 to ensure parameter $\alpha$ lying between zero and unity. The synchronization error in log-scale shows fluctuations of the order of $10^{-2}$ in Figure 9, similar to the case of sinusoidal switching. The master and slave states evolves in a same way after the initial synchronization period of around 2 sec. After that even the high rate of switching in the key parameter could not desynchronize the master and slave states (Figure 10) due to strong nonlinear structure of the global controller. The phase space diagram of the synchronization performance in this case has been shown in Figure 11.

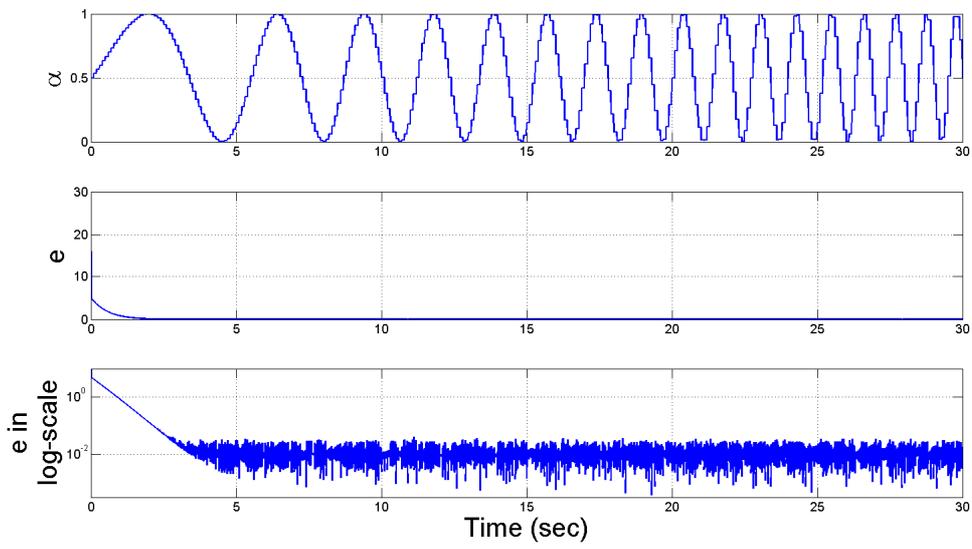

**Figure 9: chirp signal switching law for $\alpha$ and the corresponding synchronization error in normal and logarithmic scale.**

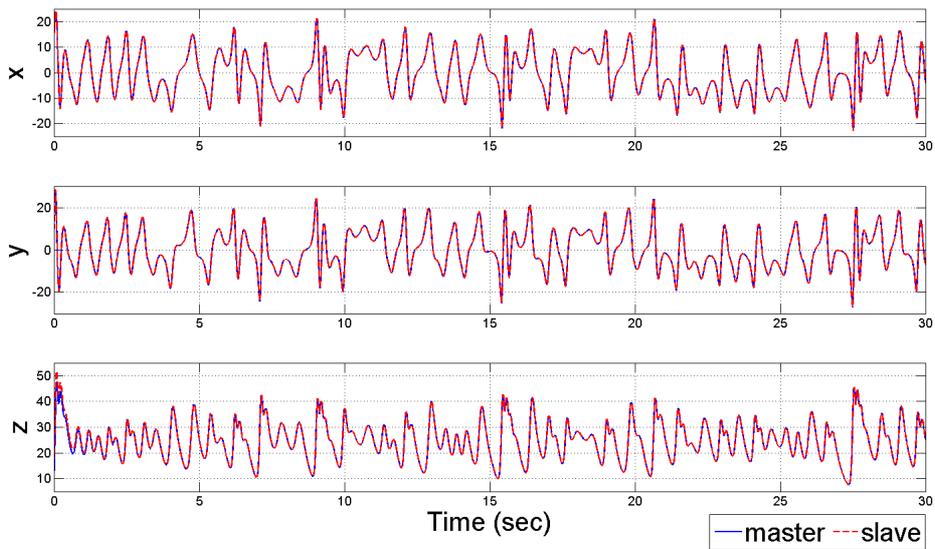

**Figure 10: Time evolution of the master-slave states for chirp signal switching in $\alpha$.**

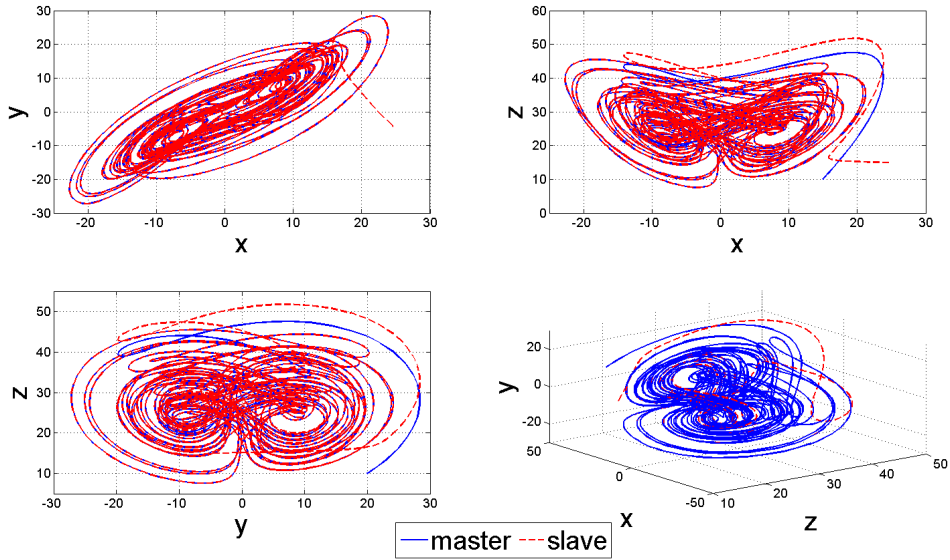

**Figure 11: Phase space representation of the master-slave states for chirp signal switching in α.**

## 4.4 Random switching in key parameter α

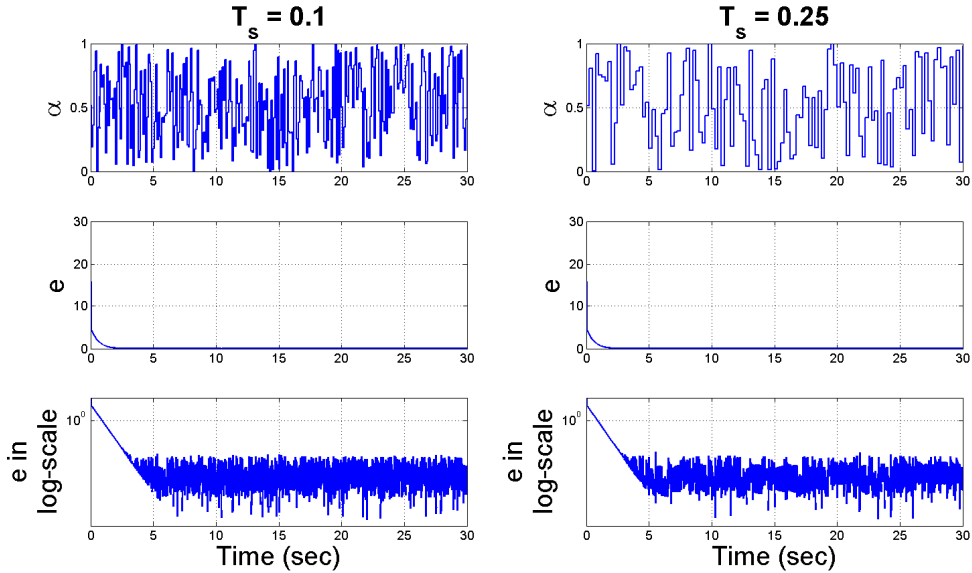

**Figure 12: Two cases of random switching law for α and the corresponding synchronization error in normal and logarithmic scale.**

Similar to the previous cases, the switching law of the key parameter $\alpha$ is chosen as a random sequence, lying between zero and unity. The random variation in $\alpha$ is sampled with a chosen sampling time of $T_s = \{0.1, 0.25\}$ sec and then passed through a ZOH block in a Matlab/Simulink based environment. The corresponding error shows acceptable synchronization performance as shown in Figure 12, even considering random switching in $\alpha$. The two cases of fast and slow switching of the key parameter in Figure 12 shows the power of the global synchronizing controller even though the characteristics of the system is continuously shuffled amongst the Lorenz, Lu and Chen family of attractors. The corresponding evolution of master and slave states in temporal and phase space domains are shown in Figure 13 and Figure 14 respectively for the two cases of fast and slow noisy sequence of the key parameter. The nature of the time evolutions of the states in Figure 13 and the nature of the phase space diagrams in Figure 14 clearly show that the dynamical behavior of the two systems under two different rate of noisy switching law are different, although in both of them the master and slave states synchronizes successfully with the same global controller.

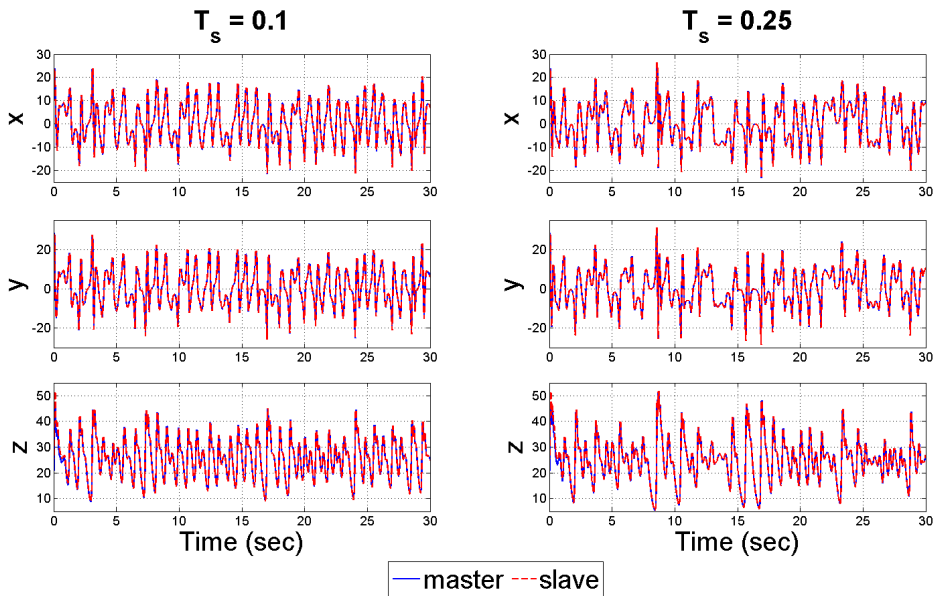

Figure 13: Time evolution of the master-slave states for fast and slow random switching in $\alpha$.

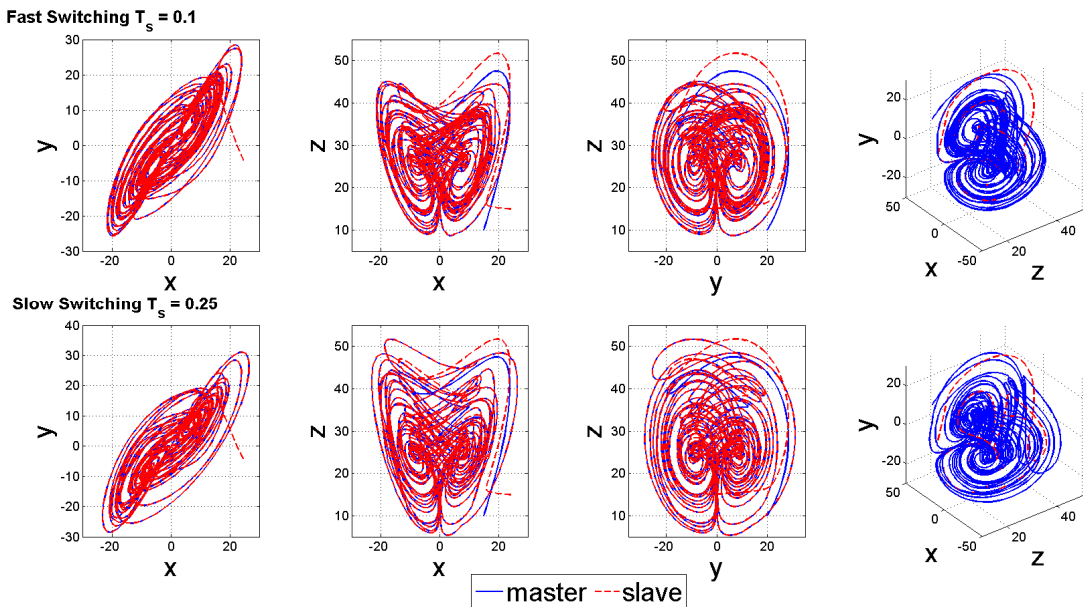

Figure 14: Phase space representation of the master-slave states for fast and slow random switching in $\alpha$.

### 4.5 Random switching in initial conditions along with random switching in key parameter α

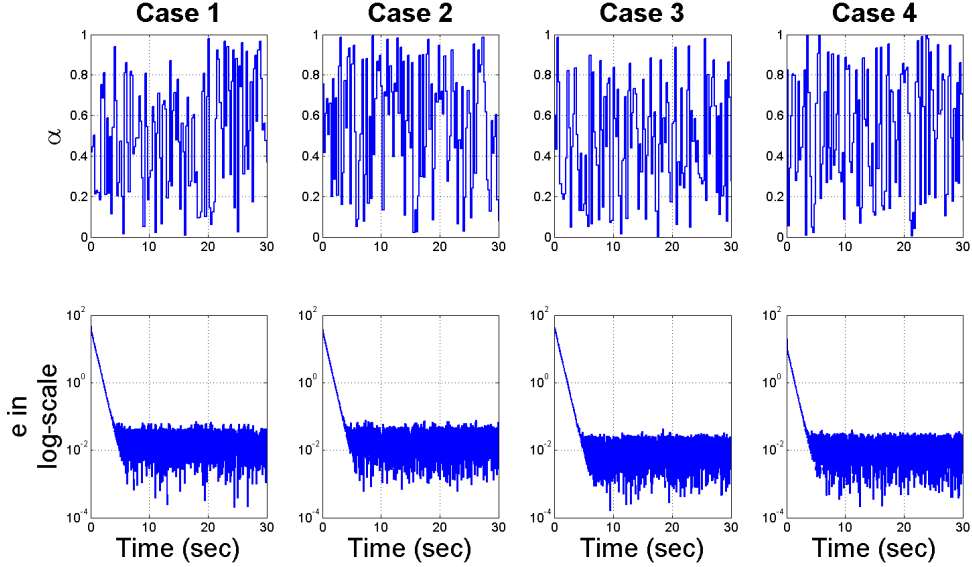

**Figure 15:** Four cases of random switching law for *α* along with random switching in initial conditions and the corresponding synchronization error in normal and logarithmic scale.

In the previous simulation examples the initial conditions for the master and slave states are chosen as $(x_{m0}, y_{m0}, z_{m0}) = (15, 20, 10)$ and $(x_{s0}, y_{s0}, z_{s0}) = (25, -5, 15)$ respectively. Here we study the synchronization performance with simultaneous random switching in key parameter *α* and the initial conditions of the state variables at the beginning of synchronization study. The goal of the study here is to observe whether the same global controller can enforce guaranteed synchronization irrespective of the choice of initial condition of the state variables. In order to simulate this condition, we explore four different cases with random initial condition for master states $(x_{m0}, y_{m0}, z_{m0}) \in [-30, 30]$ and that for the slave states $(x_{s0}, y_{s0}, z_{s0}) \in [-30, 30]$. In addition to the random initial value of master and slave states, the random sequence for the key parameter *α* has also been shuffled in different combinations (with a sample and hold of $T_s = 0.25$ sec) and we here report four representative examples (Case1 - Case4) as shown in Figure 15. Even with the random change in initial condition and the key parameter, the states of the master and slave chaotic systems are synchronized very fast as shown in Figure 16, although the dynamical behavior of the state trajectories looks different due to switching in *α*. The corresponding phase space diagrams in Figure 17 explains that the same global controller is also capable of enforcing synchronization of the master and slave states although the state trajectories could start from anywhere in the state space. The synchronization from any arbitrary initial guess of the master and slave states are guaranteed in these cases, even in the presence of a continuous random disturbance in the form of changing the underlying characteristics of the master and slave systems. From the above simulation studies it can be inferred that the proposed nonlinear controller is capable of enforcing guaranteed synchronization among the master-slave states subjected to various different switching conditions and thus can be termed as a truly global controller for these typical class of switched unified chaotic systems.

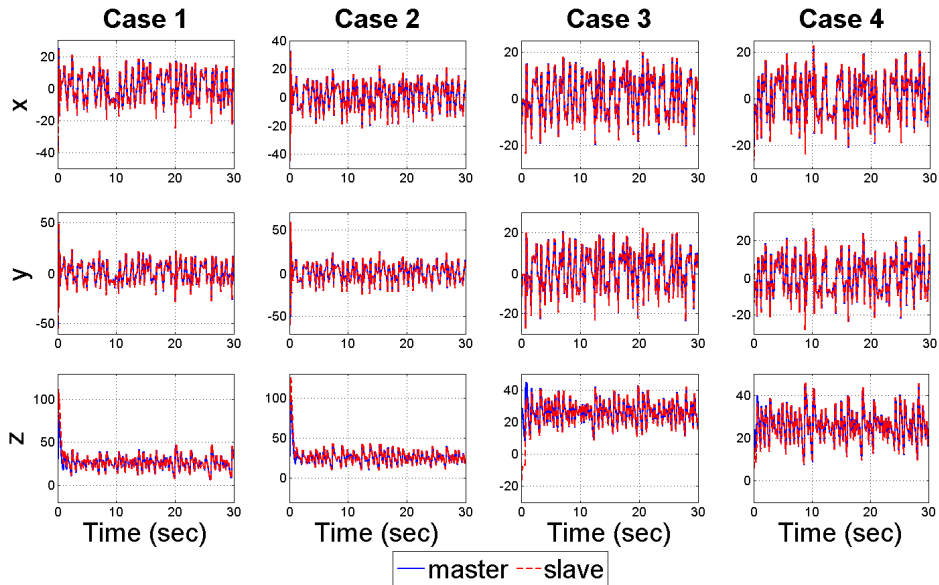

**Figure 16: Time evolution of the master-slave states for sudden switching in $\alpha$ along with random switching in the initial conditions.**

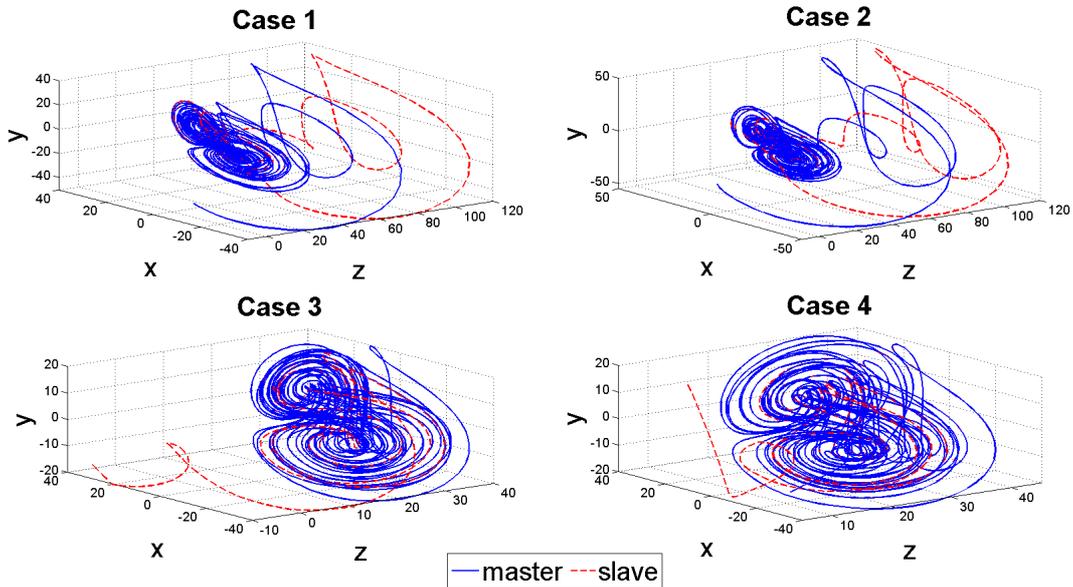

**Figure 17: Phase space representation of the master-slave states for sudden switching in $\alpha$ along with random switching in the initial conditions.**

## 4.6   Effect of sudden on-off switching of the synchronizing controller

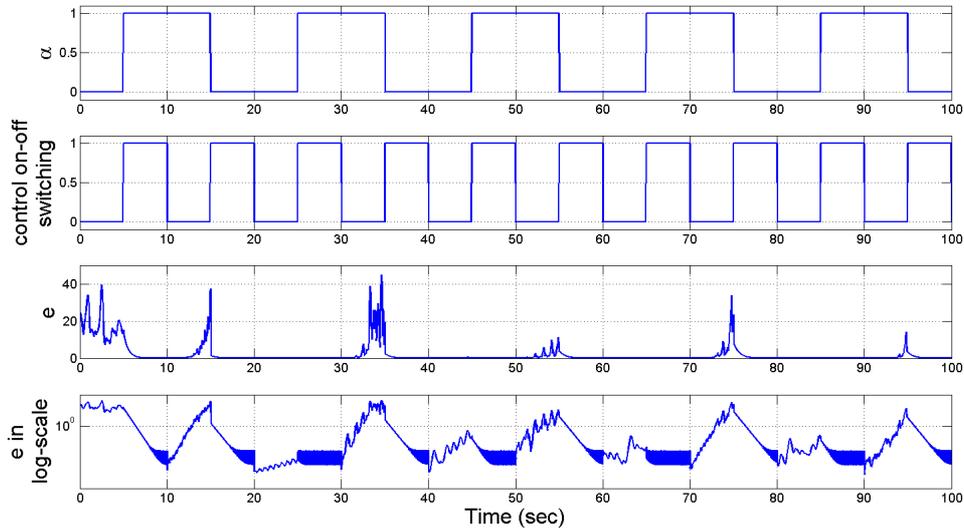

**Figure 18: Sudden step switching in α, controller on-off time sequence and the corresponding synchronization error in normal and logarithmic scale.**

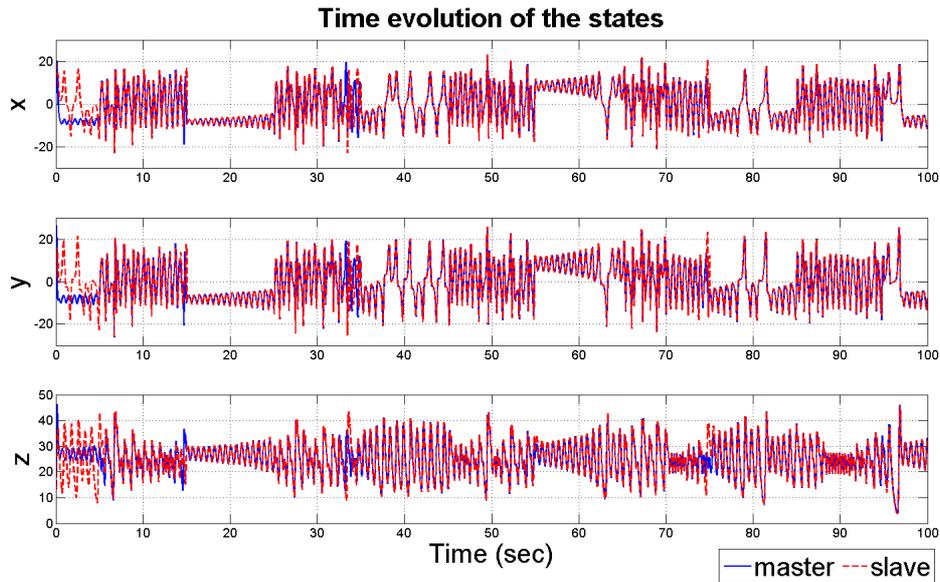

**Figure 19: Time evolution of the master-slave states for sudden on-off switching of the synchronizing controller.**

It is shown in the above simulation examples that when the controller is in action, the slave states are strongly enforced to follow the master states. Even large and rapid switching disturbances, for example, varying the key parameter $\alpha$ and changing the nature of the chaotic system, changing the initial conditions and $\alpha$ cannot desynchronize the state trajectories, exhibiting the power of the synchronizing global controller. Here we show the effect on the achievable synchronization performance if the controller is switched off for a finite duration. Intuitively it can be said that when the control action is turned off, the states of the master and slave will start to desynchronize. It may be a representative case to simulate power loss in the active control circuit or deliberate attempt to reduce the requirement of external control action for chaos synchronization. The goal of this study is to see if the states can be again synchronized without the need of designing the controller once again, after it was turned off for a while.

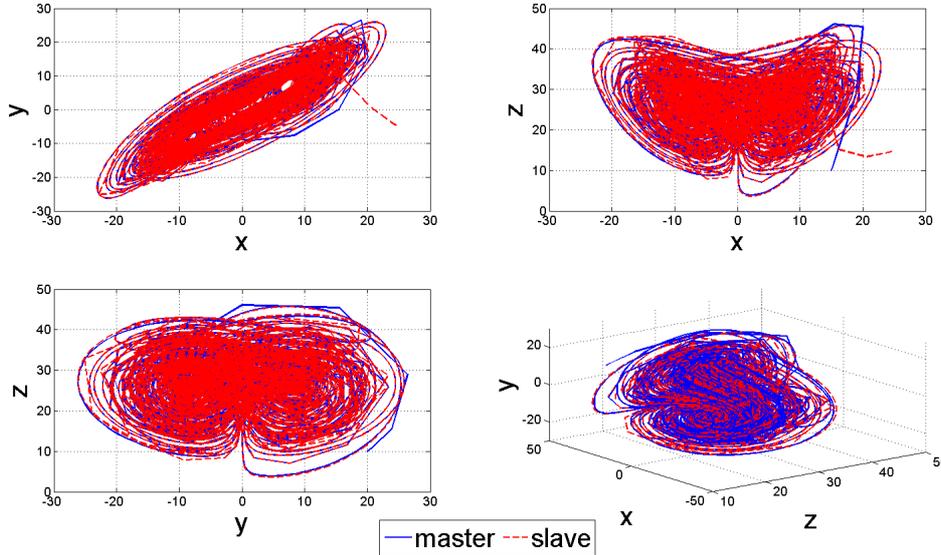

**Figure 20: Phase space representation of the master-slave states for sudden on-off switching of the synchronizing controller.**

In order to simulate this specific scenario, the key parameter $\alpha$ of the unified chaotic system is switched between zero to unity as a square wave. The duty cycle of the switching law of $\alpha$ is 50%, implying that the unified chaotic system is simulated with $\alpha = 0$ and $\alpha = 1$ for equal time scale. It is observed from Figure 18 that the period of the square wave for varying $\alpha$ is 20 sec with a phase delay of 5 sec. In order to show the desynchronization effects in the absence of the global controller, its on-off switching sequence is also represented by a square wave having a double frequency with respect to the previous case of switching in key parameter i.e. having a period of 10 sec, a phase delay of 5 sec and 50% duty cycle, as shown in the second diagram of Figure 18. From the timing diagrams in Figure 18, it is worth noticing that whenever a switching in the key parameter $\alpha$ occurs, the synchronizing controller is turned on and after every 10 sec when the error dies down, it is switched off by letting the desynchronization process start gradually. After a regular interval of 5 sec, it is again switched on and rest of the following steps occur in the same sequence as before. From Figure 18 it is observed that the synchronization error increases during the finite off-period of the controller and then rapidly falls once the controller is again switched on. This shows that even if there is a loss of power in the synchronizing control circuit, the controller need not to be tuned again even if the states of the master and slaves have been desynchronized widely during the inactive session of the controller. As a result quick and efficient synchronization can be achieve any time and the proposed controller can be considered as a plug and play device for any arbitrary combination of the switching conditions. This controller on-off strategy can also be used for efficient power management of the synchronizing control circuit when turning the controller on for the whole time period may not be a necessity and active control can only be turned on once the master and slave states have deviated to a large extent. The temporal evolution of the master and slave states are depicted in Figure 19 and the deviation or desynchronization of the two states is visually perceptible during the off period of the controller unlike the previous simulations. In the corresponding phase space diagrams in Figure 20, more rough trajectories are observed for both the master and slave during the sudden controller on-off switching. But still the global controller is capable of faithfully ensuring synchronization between the master and the slave systems.

### *4.7 Contributions and few discussions*

The main advantages of the proposed global nonlinear controller design technique for synchronization of Unified chaotic system family are as follows:

a) The synchronization is maintained under arbitrary switching between different subsystems of the Unified chaotic system family of attractors, where the master and slave systems are identical but their dynamical behavior is switched from one family to the other amongst Lorenz, Lu and Chen family. The switching laws need not be known a-priori for the controller design and can be varied as per the designer's wish.

b) The synchronizing controller design is based on the formulation of BMIs which is then transformed to few set of LMIs which are easily available using standard LMI solvers [12], [13]. Also the controller gains need to be calculated only once since the nonlinear structure is already fixed in (6) and (9). The controller as shown in the illustrative examples is capable of carrying out the following tasks:

   i. Synchronize two identical unified chaotic systems without switching.

   ii. Synchronize two identical unified chaotic systems with arbitrary switching in the key parameter $\alpha$.

   iii. Synchronize both switched and non-switched variants for arbitrary initial conditions of the master and slave unlike [38], [39].

   iv. Synchronize both switched and non-switched variants even if the controller is turned off for a finite time interval. The synchronization is established through numerical examples of periodic on-off switching of the controller having a duty cycle of 50%.

c) A single controller is able to synchronize all the above cases in (b). Hence the controller can be implemented in hardware after tuning it by solving the LMIs and ported to different chaotic systems according to the requirements of the user and the parameters need not be retuned or the numerical simulation need not be redone for the three different families of the Unified chaotic system. Thus the proposed controller can be viewed as a global controller for a wide class of switched chaotic systems having a generic structure (1) with arbitrary switching in the key parameter $\alpha$.

d) Most importantly, we here explore a new class of chaotic systems known as the 'switched chaotic systems'. In general that knowing the exact structure of the chaotic system, for a specific choice of the initial conditions the state trajectories are possible to predict in a deterministically. On contrary, here we show that even starting with the same system structure and initial conditions of the states, if the chaotic system jumps to a different family, the deterministic prediction of the states cannot be done with only the knowledge of the system structure and initial conditions. In such cases, the overall dynamics is not only guided by the chaotic system structure, parameters and initial conditions but also the switching law which can be considered as an additional key for specific applications of chaos synchronization like secure communication etc. This may help to get an extra degree of freedom for deterministic prediction of such new class of 'switched chaotic systems' using the knowledge of the switching law, along with the conventional measures like the structure and parameters of the chaotic system and the initial conditions of the three master and slave states.

## 5. Conclusion

A nonlinear state feedback law and a LMI based global controller design technique is proposed here for synchronization of a new class of switched chaotic systems. The paper proposes a single controller which can synchronize two identical master-slave chaotic systems which arbitrarily switches between different variants of the family of Unified chaotic systems *viz.* Lorenz, Lu and Chen family. The same controller is also able to synchronize between identical unified chaotic systems without any switching. Illustrative examples with credible numerical simulation show the effectiveness of the proposed approach for sudden step switching, sinusoidal, chirp signal and fast and slow random switching in chaotic system parameters along with arbitrary initial conditions of the master and slave states. The controller can also enforce fast synchronization even if it is periodically turned on and off, thus acting as a truly global controller for a wide class of switched chaotic systems. Future scope of research can be directed towards synchronization and control of switched hyper-chaotic and fractional order chaotic systems.